\documentclass[12pt]{amsart}
\usepackage{amsmath}
\usepackage{enumerate}
\usepackage{amssymb}
\usepackage{amsbsy}
\usepackage{amsfonts}

\theoremstyle{plain}
\newtheorem{thm}{Theorem}[section]
\newtheorem{prop}[thm]{Proposition}

\theoremstyle{definition}

\newtheorem{rmk}[thm]{Remark}
\numberwithin{equation}{section}
\newcommand{\sm}{\left(\begin{smallmatrix}}
\newcommand{\esm}{\end{smallmatrix}\right)}

\newfont{\FieldFont}{msbm10 scaled\magstep1}

\newcommand{\pf}{\noindent\bf Proof }

\newtheorem{theorem}{Theorem}[section]

 \newtheorem{example}[theorem]{Example}

\numberwithin{equation}{section}

\newcommand{\blanks}[1][4em]{\underline{\makebox[#1]{}}}

\usepackage{color}

\definecolor{blue}{rgb}{0,0,1}
\definecolor{red}{rgb}{1,0,0}
\definecolor{green}{rgb}{0,.6,.2}
\definecolor{purple}{rgb}{1,0,1}

\long\def\red#1\endred{{\color{red}#1}}
\long\def\blue#1\endblue{{\color{blue}#1}}
\long\def\purple#1\endpurple{{\color{purple}#1}}
\long\def\green#1\endgreen{{\color{green}#1}}

\begin{document}

\title{Periods of Hilbert Modular forms,  Kronecker series  and Cohomology }
\author{YoungJu Choie  }

\address{YoungJu Choie \endgraf
Department of Mathematics\endgraf
Pohang University of Science and Technology (POSTECH) \endgraf
Pohang,
Republic of Korea
}
\email{yjc@postech.ac.kr}
\today

 \thanks{This work was partially supported by
 NRF 2018R1A4A1023590 and NRF 2017R1A2B2001807}

 \subjclass[2000]{ 11F41, 11F50, 11F60, 11F67}
 \keywords{Hilbert modular form, parabolic cohomology, period polynomial,  Kronecker series }
\maketitle
\begin{abstract}
 Generalizing a result of~\cite{Z1991, CPZ} about elliptic modular forms, we give a closed formula for the sum of all Hilbert Hecke eigenforms over a totally real number field with  strict class number $1$, multiplied by their period polynomials, as a single product of the Kronecker series.  
\end{abstract}

\section{\bf{ Introduction }}

Based on Bol's result \cite{Bol} Eichler  initiated a theory of the periods of integrals so that an
  automorphic form of the first or second kind leads to a cohomology class in the mapping of a Fuchsian group into a polynomial module and the (converse) correspondence of each such cohomology class leads to an automorphic form in one complex variable \cite{Eichler}.
Shimura extended this theory    by showing that   the structure of an abelian variety in certain cases can   be also   given to the periods of such integrals and  showed   critical values of the $L$-functions attached 
to elliptic modular forms can be computed explicitly using the cohomology group \cite{Shimura1959}.
This method was developed by Manin \cite{M} who  proved an algebraic theorem for the periods of elliptic cusp forms  for the full modular group and studied $p$-adic properties of the algebraic factors in $L$-functions.
Kohnen-Zagier \cite{KZ}  further extended this theory to elliptic modular forms including Eisenstein series \cite{KZ} and studied forms whose  period polynomials  have arithmetically interesting  rational structure relating to Bernoulli numbers,   binary quadratic forms,   zeta-functions of real quadratic fields,   modular forms of half-integral weight  and   Hilbert modular forms.
Hence,  period polynomials, which
  allow  us to compute the critical values of  $L$-function of modular forms at once,
give a rich source of relations between modular forms and arithmetic.

\medskip

The period polynomial of an elliptic cusp form $f(\tau)=\sum_{\ell\geq 1}a_f(\ell)q^{\ell}\, \,  (\tau \in \mathbb{H}=$ upper half plane, $q=e^{2\pi i \tau})$ of weight $k$ on $SL_2(\mathbb{Z})$ is the polynomial of degree $k-2$ defined by
\begin{eqnarray}\label{pe}
r_f(X)=\int_{0}^{i\infty}f(\tau)(\tau-X)^{k-2}d\tau
\end{eqnarray}

or equivalently  by
$$r_f(X)=-\sum_{n=0}^{k-2}\frac{(k-2)!}{(k-2-n)!} \frac{L(f,n+1)}{(2\pi i)^{n+1}} X^{k-2-n}, $$
where $L(f,s)=\sum_{n\geq 1}\frac{a_f(n)}{n^s} \,  (Re(s)\gg 0) .$
The maps $f\rightarrow r_f^{ev}$ and $f\rightarrow r_f^{od}$ assigning to $f$ the even and odd parts of $r_f$ are both injective with known images  from the  Eichler-Shimura-Manin theory.

\medskip
When $f$ is a Hecke eigenform then  one has the   two-variable polynomial
\begin{eqnarray*}
r_f(X,Y):=\frac{r_f^{ev}(X)r_f^{od}(Y) + r_f^{od}(X)r_f^{ev}(Y)}{(2i)^{k-3} <f, \, f>} \in \mathbb{Q }_f [X,Y]
\end{eqnarray*}
where $ \mathbb{Q}_f$ is the field generated by Fourier coefficients of $f$ over $\mathbb{Q}.$
\medskip

 Zagier \cite{Z1991} found the following attractive formula :

\begin{eqnarray}\label{Za1991}
&&    \frac{(XY-1)(X+Y)}{X^2Y^2}T^{-2}
 +  \sum_{k\geq 2}
\sum_{  f \in \mathcal{B}_k
 }
r_f(X,Y) f(\tau) \frac{T^{k-2}}{(k-2)!}
\\
&=&
F_{\tau}(T, -XYT) F_{\tau}(XT, YT),  \, F_{\tau}(u,v)=
\frac{
\theta^{'}(0)
\theta(u+v)}
{\theta(u)\theta(v) }\nonumber
 \end{eqnarray}

\noindent
where
 $ \theta(u)= \sum_{n\in\mathbb{Z}}(-1)^nq^{\frac{1}{2}(n+\frac{1}{2})^2}e^{(n+\frac{1}{2})u} $ is  
a Jacobi theta function   and
$\mathcal{B}_k$ is a set of all Hecke eigenforms of weight $k$ on $SL_2(\mathbb{Z}).$

The identity by Zagier (\ref{Za1991})  relates a generating function, which contains  all Hecke eigenforms  together with  all   critical values,  to   the Jacobi  form $F_{\tau}(u,v). $
Such expansions with respect to the variable  $T$   give  an algorithm to compute Hecke eigenforms (see \cite{Z1991}  for  more details).
\\

It took almost 30 years to see that such an identity (\ref{Za1991}) is not accidental but
 also exists for a general group $\Gamma_0(N)$ (see \cite{CPZ}).
Now it seems natural to ask if one can get such a relation for general automorphic forms.
In this paper  we  attempt to get such a formula, namely,  an identity between
a generating function of periods of    Hilbert modular forms   over totally real number fields with   strict class number one  and Jacobi forms (see Theorem \ref{main}).   The function
$F_{\tau}(u,v) $ in (\ref{Za1991}) was introduced by Kronecker  in a more general form (see page 70 in \cite{Weil}) and  several properties of that have been explored by Zagier \cite{Z1991}.  The essential property of   $F_{\tau}(u,v) $ is that  it
can be identified as a sum of derivatives of Eisenstein series (called  the  "Kuznetsov lifting") and  using this fact we are  able to extend  Zaiger' identity to  that for a totally real  number field.   The main result of this paper shows the first connection between the Kronecker series and the critical $L$- values of     Hilbert modular forms over the totally real number fields.   It also gives a systematic way  to compute  Hilbert Hecke eigenforms and  the special values of $L$ -functions  by taking the expansions of the Kronecker series.   
\medskip

 This paper is organized as follows: in  section $2$ we state (Main) Theorem  \ref{main} after introducing 
necessary notations.  
In section $3$ the analog  of Eisenstein-Kronecker series over a totally real    number field and  the rationality of period polynomials of Hilbert modular forms are  discussed.  Section $4$ gives   detailed proof of Main Theorem.
Finally, we give a comment on a connection between   parabolic cohomology  and a period theory of Hilbert modular forms.  Also, we discuss  a  possible application of the Kronecker series to evaluate the  special $L$-values of  a  general automorphic form as a conclusion.

\medskip

\noindent {\bf{Acknowledgement}}  I would like to thank the
referees for numerous helpful comments and suggestions which greatly
improved the exposition of this paper.

\medskip
\section{\bf{Notations and statement of Main Theorem }}
 \subsection{Notation}

\begin{itemize}
\item
  $\mathbb{F} :$
a totally real number field of degree $t$  with discriminant  $D$ and  strict class number $1$
\item    $ \mathcal{O} :$   the ring of
integers of $\mathbb{F} $ containing a unit of negative norm
\item  $\mathcal{O}^*: $ the group of  units  of  $ \mathcal{O}$
 \item $\mathcal{O}^{*,+} : $  the group of totally positive units of $\mathcal{O} $
\item  $U^+=\{\epsilon^2\, : \, \epsilon\in  \mathcal{O}^{*,+}\} $
\item $\alpha \succ 0 : \, $ a totally positive element
\item   $\alpha_1 \rightarrow \alpha_2,   \cdots, \alpha_t$ for the conjugation
\item   $\mathcal{N}(\alpha)=\prod_{j=1}^{t} \alpha_j$ the norm
\item   $tr(\alpha)=\sum_{j=1}^t \alpha_j $ the trace
\item  $\mathfrak{D} : $ the different of $\mathbb{F}$ 
 \item  $\zeta_{\mathbb{F}}(s) = 
\sum_{ {c}\subset \mathcal{O} } \frac{1}{\mathcal{N}( {c})^s}, $
where $s\in \mathbb{C} $ and the  integral ideal $c $

\item  $|\mathbf{r}|= \sum_{i=1}^t r_i, \, \mathbf{r}+\mathbf{r}'=(r_1+r_1', \cdots, r_t+r_t'),\,$  
$  \Gamma(\mathbf{r}+\mathbf{1})= \mathbf{r}!=r_1! \cdots r_t!, \,  \sm \mathbf{r} \\ \mathbf{r}'\esm = \sm r_1  \\  r_1' \esm  \sm r_2  \\  r_2' \esm
\cdots \sm r_t  \\  r_t' \esm$ for $\mathbf{r}=(r_1,  r_2, \cdots, r_t),  \mathbf{r}'=(r_1', r_2', \cdots, r_t')  \in \mathbb{Z}_{\geq 0}^t   $ 

 \item $z^{\mathbf{r}} = {z_1}^{r_1}{z_2}^{r_2}\cdots{z_t}^{r_t},\, tr(\mathbf{m}z) = \sum_{j=1}^t m_j z_j, \, \mathcal{N}(z) = \prod_{j=1}^t z_i$ for $z = (z_1, z_2, \cdots, z_t)\in \mathbb{C}^t $ and ${\mathbf{m}} \in  \mathbb{F}$
\item  $\mathbb{H}^t : $ the $t$-copies of complex upper half plane $\mathbb{H}$
\item $\tau=(\tau_1, \cdots, \tau_t)=x+\sqrt{-1}y \in \mathbb{H}^t, x=(x_1, \cdots, x_t) \in \mathbb{R}^t,  y = (y_1, \cdots, y_t) \in  (\mathbb{R}^+)^t , \, q = \prod_{j=1}^t q_j, \, q_j=e^{2\pi i \tau_j}, \, 1\leq j\leq t.$
\item $ \sigma=(\sigma_1, \cdots,  \sigma_t) \in  \Gamma=SL_2(\mathcal{O})^t : $ an element in the Hilbert modular group
\item
The action of the group $\Gamma$, which is embedded into
 $SL_2(\mathbb{R}) \times \cdots \times SL_2(\mathbb{R}),$ on $\mathbb{H}^t$ is
given by  linear fractional  transformations
$$\sm a&b\\c&d\esm \tau=\frac{a\tau+b}{c\tau+d}=\bigl(\frac{a_1\tau_1+b_1}{c_1\tau_1+d_1}, \cdots,  \frac{a_t\tau_t+b_t}{c_t\tau_t+d_t}\bigr), \tau=(\tau_1, \cdots, \tau_t)\in \mathbb{H}^t $$
\item
For a holomorphic function $\chi$ on $\mathbb{H}^t,$
$$\chi^{(  \ell )}(\tau)=\frac{\partial^{| \ell |}}
{\partial{\tau }^{\ell } } \chi(\tau)
:=\frac{\partial^{| \ell |}}
{\partial{\tau_1}^{\ell_1}\cdots\partial{\tau_t}^{\ell_t}}\chi(\tau),\, \, \,  \forall \ell=(\ell_1, \cdots, \ell_t)\in \mathbb{Z}^t_{\geq 0} 
$$
\item  $ \mathbb{D}^{ \ell }\bigl(\chi(\tau)\bigr):=
  \chi^{(\ell)}(\tau)$

\item $S_{ \mathbf{k}}\subset M_{\mathbf{k}} :$  the space of Hilbert cusp form $\subset $ the space of Hilbert modular form on $\Gamma$ with  a parallel weight  $\mathbf{k}=(k,\cdots, k), $ even $k\geq 2.$

 \item $\mathcal{B}^{0}_{k}\subset \mathcal{B}_{k}: $ a basis, consisting of all normalized Hecke eigenforms, of $S_{ \mathbf{k}}\subset M_{ \mathbf{k}} ,$ respectively

\item $\mathbb{Q}_f : $  the field spanned by Fourier coefficients of $f$ over $\mathbb{Q}$
\item
$
\bigl< f\, ,\, g \bigr> := \int_{\Gamma \backslash \mathbb{H}^t} f(\tau)\overline{g(\tau)}\frac{dx\, dy}{y^2} , \,
$
the Petersson inner product for $f \in S_{\mathbf{k}}, g\in M_{\mathbf{k}}$
\item For a function $f$ on $\mathbb{H}^t$ and $\mathbf{\ell}=(\ell_1, \cdots, \ell_t) \in \mathbb{Z}^t,$
\begin{eqnarray*} (f|_{\ell}\sigma)(z):=(cz+d)^{-\mathbf{\ell}}f(\frac{az+b}{cz+d}),   \sigma =\sm a & b\\ c& d\esm \in \Gamma
\end{eqnarray*}
\end{itemize}

\subsection{{\bf{Statement of Main Theorem}}} 
 
Take   a cusp form $ f(\tau ) = \sum_ {\mathfrak{D}^{-1} \ni \nu \succ  0}a_f(\nu) e^{2\pi i tr(\nu \tau)} $ in $ S_{\mathbf{k}}$ and
  consider the complete $L$-function of $f:$  for $s\in \mathbb{C}, $
$$\Lambda(f, s):=\int_{\mathbb{R}_+^{t} /  U^+  } f(iy )
  y^{s-1} \, dy
=D^{s}(2\pi)^{-ts }
\Gamma(s )^t  L(f, s ) ,$$
where $\, L(f,s)= \sum_{ \mathfrak{D}^{-1}/  U^+ \ni \nu \succ 0  }
\frac{a_f(\nu)}{\mathcal{N}(\nu)^s} \, \, \, ( Re (s) \gg 0 ).$   \, 
It is well-known that $\Lambda (f, s)$ has an analytic continuation and functional equation \cite{123, Geer}
\begin{eqnarray*}\label{functional-1}
\Lambda (f, s)=  (-1)^{\frac{tk}{2}} \Lambda (f, k-s).
\end{eqnarray*}

Consider the following polynomials in $X=(X_1,\cdots, X_t)$, called the {\bf{ even (odd) period polynomial}} associated to $f :$
\begin{eqnarray*}\label{poly}
R_f^{ev} (X)&:= &\sum_{\tiny{\begin{array}{cc} 0\leq n\leq k-2 \\ n\equiv 0 \pmod{2} \end{array}}}
\frac{\Gamma({k}- {1})^t}{\Gamma(n+1)^t\Gamma(k-n-1)^t}  R_{ {k-2-n}}(f) \mathcal{N}(X)^{  {n}},\nonumber\\
R_f^{od} (X)&:= &\sum_{\tiny{\begin{array}{cc} 0< n < k-2 \\ n\equiv 1 \pmod{2} \end{array}}}
\frac{(-1)^{nt}\Gamma({k}- {1})^t}{\Gamma(n+1)^t\Gamma(k-n-1)^t}  
  R_{ {k-2-n}}(f)\mathcal{N}(X)^{  {n}}, \nonumber\\
 && \nonumber \\
& &  R_f(X):=  (-1)^t \bigl(R_f^{ev} (X) + R_f^{od} (X)\bigr), \, \,
\end{eqnarray*}
 where
\begin{eqnarray*}\label{value}
&& R_{ {n}}(f): =\int_{{\mathbb{R}_{+}}^t /  U^+ } f(\tau)\tau^{n}d\tau  =  i^{t(n+1)}\Lambda(f,n+1).
\end{eqnarray*}
 Using the functional equation of $\Lambda(f,s)$  we get
\begin{eqnarray*}
R_{{k-2-n}}(f)
=(-1)^{t(n+1)}R_{ { n}}(f)
\end{eqnarray*}
and  so we get
$$ \mathcal{N}(X)^{k-2}R_f(-\frac{1}{X }) =(-1)^{t}R_f(X).$$
\medskip
Let $f$ be a primitive (Hilbert) Hecke eigenform and
consider the  polynomial of  the $2t$-variables in $X=(X_1, \cdots, X_t)$ and $ Y=(Y_1, \cdots, Y_t)$
\begin{eqnarray*}
R_f(X,Y):=(-1)^t  \frac{R^{ev}_f(X) R^{od}_f(Y)+R^{ev}_f(Y) R^{od}_f(X)}
{  D^{k-\frac{1}{2}} \,  (2i)^{t(k-3)}  <f, \, f>} \in \mathbb{C}[X,Y] .
\end{eqnarray*}
It  transforms under $\sigma\in Gal(\mathbb{C} / \mathbb{Q})$ by $R_{\sigma(f)}=\sigma(R_f) $ so that $R_f(X,Y)$ has   coefficients in the number field $\mathbb{Q}_f$ generated by the Fourier coefficients of $f.$ 
Summing over  the basis $\mathcal{B}_k^0,$  consisting of all normalized Hecke eigenforms of $S_{ \mathbf{k}},$
the following function
\begin{eqnarray}\label{cc}
C_k^{cusp} (X,Y;\tau ) := \sum_{
f\in \mathcal{B}_k^0
}R_f(X,Y)f(\tau)
\end{eqnarray}
is in $\mathbb{Q}[[q]][X, Y]$ for each even integer $k\geq 2.$
Further we extend the definition of  $R_f(X,Y)$  (see  section \ref{non})  to non-cusp forms and include the Eisenstein series in the sum (\ref{cc}). Then we define 
\begin{eqnarray}\label{kge}
C_k (X,Y;\tau ) := \sum_{f\in \mathcal{B}_k}
R_f(X,Y)f(\tau). \end{eqnarray}

\begin{example} Take $\mathbb{F}=\mathbb{Q}(\sqrt{5}).$
\begin{enumerate}
\item 
\begin{eqnarray*}
&&C_2(X, Y;\tau)\\
&&= \sum_{
f\in \mathcal{B}_k
}R_f(X,Y)f(\tau)=2^4 \cdot3\cdot 5
\frac{ \bigl(\mathcal{N}(X)+\mathcal{N}(Y) \bigr) \bigl(\mathcal{N}(XY)+1 \bigr) }
  {\mathcal{N}(XY)}
G_{\mathbb{F}, 2}(\tau),
\end{eqnarray*}
\noindent with  the normalized Eisenstein series of weight $(k,k)$ on $\Gamma $ given by
\begin{eqnarray*}
G_{\mathbb{F},k}(\tau)= \frac{\zeta_{\mathbb{F}}(1-k)}{2^2}
+\sum_{  \mathcal{D}^{-1} \ni  \nu\succ 0}
\sigma_{k-1}(\nu \mathcal{D})e^{2\pi i tr(\nu\tau)},\sigma_{r}(\frak{n})=\sum_{\frak{c}|\frak{n}} \mathcal{N}(\frak{c})^{r}.
\end{eqnarray*}

\item Let $\mathbb{F}=\mathbb{Q}(\sqrt{5})$ and take a unique  cusp form $f$ of weight $8 $  on $\Gamma.$
Using the example  in \cite{BRW}   we get 
 $$ \frac{R_f^{ev}(X)R_f^{od}(Y)}{5^{\frac{15}{2}}(2 i)^{10} <f,f>} =
c(1+\frac{361}{2^2\cdot 5}X^2+\frac{361}{2^2\cdot 5}X^4+X^6)
(Y+\frac{2}{3}Y^3+Y^5) $$ up to rational constant multiple $c.$
\end{enumerate}  
 \end{example}
 
\medskip

Combining all these functions into a single generating function to define
\begin{eqnarray}\label{cccc}
&&  C(X,Y;\tau;T)\\
&& :=\frac{(\mathcal{N}(X)+\mathcal{N}(Y))(\mathcal{N}(X Y)+(-1)^t) }{  \mathcal{N}(X Y T)^2}+\sum_{k\geq 2}C_k(X,Y;\tau )\frac{\mathcal{N}(T)^{ {k-2}}}{\Gamma(k-1)^t}.\nonumber
\end{eqnarray}

\medskip
On the other hand,
consider
\begin{eqnarray}\label{theta1}
\, \, F_{\tau}(u,v) :=
  (-2)^t \sum_{k\geq 0}\widetilde{{G}_{\mathbb{F}, {k}}}(\tau, \frac{ u v  }{ 2\pi i } )( \mathcal{N}(u)^{  {k-1}}+\mathcal{N}(v)^{  {k-1}}), u, v \in \mathbb{C}^t,
\end{eqnarray}
 where
\begin{eqnarray}\label{Ku}
&&\\
&&
\widetilde{{G}_{\mathbb{F}, {k}}}(\tau, \lambda) :=
\biggl\{
\begin{array}{cc}
\sum_{ 
 \ell=(\ell_1,\cdots, \ell_t)\in  \mathbb{Z}^t_{\geq 0} } 
\frac{  {\lambda }^{ \ell } }{ {\ell }!   ( {\ell }+\mathbf{ {k}- {1}})!  }
\mathbb{D}^{  \ell  }\bigl({G}_{\mathbb{F}, {k}} (\tau) \bigr) & \mbox{ if $k\geq 2$}\nonumber\\
\frac{1}{2^t  } & \mbox{if $k=0$}
\end{array}
\biggr\}, \lambda \in \mathbb{C}^t 
\end{eqnarray}
  and   a  normalized Hilbert Eisenstein series $G_{\mathbb{F},  {k}}(\tau)$ (p 20  in  \cite{Geer}) defined by
 \begin{eqnarray*}\label{Eisen}E_{\mathbb{F}, {k}}(\tau) &:=& \frac{D^{ \frac{1}{2}-k }(2\pi i)^{tk} }{\Gamma(k)^t} \bigl(\frac{1}{2^t}\zeta_{\mathbb{F}}(1-k)+\sum_{\tiny{ \begin{array}{cc}\nu \in \mathfrak{D}^{-1}\\ \nu \succ  0 \end{array}}}\sigma_{k-1}(\nu \mathfrak{D})e^{2\pi i tr(\nu \tau)}\bigr)\\&:=&\frac{D^{ \frac{1}{2}-k }(2\pi i)^{tk} }{\Gamma(k)^t}  G_{\mathbb{F},  {k}}(\tau),   \, \sigma_{r}(\frak{n})=\sum_{\frak{c}|\frak{n}} \mathcal{N}(\frak{c})^{r}.\\
\end{eqnarray*}

\noindent Now we state the  main result of  this paper :
\begin{theorem}\label{main}{{\bf (Main Theorem)}}
Let   $ C(X, Y;\tau; T)$ be the generating function of the periods of Hilbert modular forms given in (\ref{cccc}). Then we have
\begin{enumerate}
\item
$ C(X, Y;\tau; T)\in \frac{1}{\mathcal{N}(XYT)^2}  {\mathbb{Q}}[X,Y][[q,T]].$   
\item 
$ 
C(X,Y;\tau;T) = F_{\tau}(T, - XYT) \, F_{\tau}(XT,YT).
$ 
\end{enumerate}
\end{theorem}
\medskip
\begin{rmk}\label{zeta}
$\widetilde{{G}_{\mathbb{F}, {k}}}(\tau, \lambda)$ in (\ref{Ku}) is
the  "Kuznetsov lifting" of  the Hilbert Eisenstein series $ G_{\mathbb{F}, \mathbf{k}}(\tau).$ Its modular transformation property is known (see Theorem 2 in \cite{CKR}) :
 for any
 $\sm a & b\\ c& d\esm \in \Gamma, k\geq 2,$
$$
\mathcal{N}(c\tau+d)^{-{k}}e^{-tr( \frac{c\lambda}{c\tau+d})} \widetilde{{G}_{\mathbb{F}, {k}}}(\frac{a\tau+b}{c\tau+d}, \frac{\lambda}{c\tau+d})=\widetilde{{G}_{\mathbb{F}, {k}}}(\tau, \lambda)$$
and    its generating function $\mathcal{F}_{\tau}(u,v):= (-2)^t \sum_{k\geq 2}\widetilde{{G}_{\mathbb{F}, {k}}}(\tau, \frac{ u v  }{ 2\pi i } )( \mathcal{N}(u)^{  {k-1}}+\mathcal{N}(v)^{  {k-1}})$   behaves as  a Jacobi-like form  \cite{EZ, CL} with a modular transformation property
$$  \mathcal{F}_{\frac{a\tau+b}{c\tau+d}}(\frac{u}{c\tau+d}, \frac{v}{c\tau+d})
=\mathcal{N}(c\tau+d)\,   e^{ tr(\frac{cuv}{  c\tau+d })}
{\mathcal{F}}_{\tau}(u,v),   \forall \sm a&b\\c&d\esm  \in  \Gamma. $$
\end{rmk}
\bigskip

\section{{\bf{ Algebraicity and Period of Hecke eigen forms }}}
\subsection{\bf{Algebraicity}}
The study of period relation for automorphic forms was started by Shimura. He showed the existence of relations up to factors in $\overline{\mathbb{Q}}^* $ in many instances  and made a general conjecture relating   periods of Hilbert modular varieties and their compact analogs, that is, the quaternionic modular varieties \cite{  Shimura-inv, Shimura-Ame}.
There is a weaker conjecture, which gives a relation
between a product of two periods, called the quadratic periods, may be interpreted, up to algebraic factors, as Petersson inner products.  This was proved by  M. Harris \cite{Ha} under a certain technical condition.
 More precisely,
for each $m, 0 \leq m \leq k-2,$
$\Lambda(f, m+1)$ is called the {\bf{critical}} values.

\begin{thm}(Theorem 4.3 in \cite{Shimura-Duke})\label{Rationality} Let $f$ be a   Hilbert Hecke eigenform   of weight  $\mathbf{k}=(k, \cdots, k)$ over a   totally real   number field $\mathbb{F}$ of degree $t$
and  $\sigma \in Gal(\overline{\mathbb{Q}}/\mathbb{Q}).$
\begin{enumerate}
\item For each $ r \in  \mathbb{Z}^t  / 2\mathbb{Z}^t$ and for $f^{\sigma}, \sigma \in Gal(\overline{\mathbb{Q}}/ \mathbb{Q}),$
there exist nonzero complex numbers
$\omega_{f}^{r} $  such that
 $ (\frac{L(f,m)}{(2\pi i)^{tm} \omega_f^{r} })^{\sigma}
=\frac{L(f^{\sigma}, m)}{(2\pi i)^{tm} \omega_{f^{\sigma}}^{r} },$  for any integer $m$
such that  $0<m<k.$

\item $\frac{L( f, m)}{(2\pi i)^{tm} \omega_f^{r} } \in \mathbb{Q}_f.$

\item If $p=(p_1, \cdots, p_t), r=(r_1, \cdots, r_t)$ with $p_i+r_i\equiv 1\pmod{2}, 1 \leq \forall i \leq t, $
we have
 $\frac{w_f^p \cdot w_f^r }{ <f\, , f> } \in \mathbb{Q}_f $
and
 $\bigl(\frac{w_f^p \cdot w_f^r }{ <f\, , f> }\bigr)^{\sigma}
=    \frac{w_{f^{\sigma} }^p \cdot w_ {f^{\sigma}}^r }{ <f^{\sigma}\, , f^{\sigma}> }.  $

\end{enumerate}
\end{thm}

\subsection{{\bf{ Period of non-cusp forms }}} \label{non}
 

Since $\mathcal{B}_{k}$ in (\ref{kge})  contains non-cusp forms   one needs to explain  "{\em{ period function}}" corresponding to a non-cusp  form $f.$
  \noindent Take a non-cusp form  \mbox{ $f(\tau ) = \sum_{ 
 0 \preceq \nu \in \mathfrak{D}^{-1} 
 } a_f(\nu) e^{2\pi i tr(\nu \tau)} $}  in $M_{\mathbf{k}}$ and consider 
$$\Lambda(f,s):=\int_{( \mathbb{R}_{+})^t /  U^+  }\bigl(f(iy ) - a_f(0)\bigr) y^{ {s}-1 } dy
=D^{s}(2\pi)^{-ts}\Gamma(s)^{t} L(f, s),   s\in \mathbb{C}. $$
 It has a meromorphic continuation to $\mathbb{C}$ and satisfies a functional equation $\Lambda(f, s) = (-1)^{\frac{tk}{2}}\Lambda(f, k-s), $
but now has simple poles of residue  as    $ - a_f(0)$ and $(-1)^{kt}a_f(0)$ , up to a constant multiple, at $s=0$ and $s=k,$ respectively.

\noindent Define
\begin{eqnarray}\label{Pe}
&&R_f(X):=\frac{
(-1)^{t }\sqrt{D} \cdot a_f(0)}{ ( {k-1} )^t}(\mathcal{N}(X)^{ {k-1}}+(-1)^{t}
 \mathcal{N}(X)^{- {1}}) \\
&+& \sum_{n=0}^{k-2} (-1)^{\frac{t(k+n-1)}{2}}
 \frac{\Gamma({k-1})^t}{ \Gamma(n+1)^t\Gamma(k-n-1)^t} \Lambda(f, k-1-n) \mathcal{N}(X)^{ n}.  \nonumber
\end{eqnarray}

\bigskip

The assumption that $\mathbb{F}$ has the strict class number $1$ implies 
 that  the space of Hilbert modular forms  is a  direct sum    (see   \cite{123}, p 12)  $$ M_{\mathbf{k}} = S_{\mathbf{k}} \oplus  < G_{\mathbb{F}, k} > $$
and, so  (\ref{kge}) becomes
$$C_{k}(X,Y;\tau) = \sum_{f \in \mathcal{B}_k^0} R_f(X,Y)f(\tau) 
+R_{G_{\mathbb{F},k}}(X,Y) G_{\mathbb{F},k}(\tau) , $$  where 
  $R_{G_{\mathbb{F},k}}(X, Y),$  defined by
\begin{eqnarray}\label{pee}
R_{G_{\mathbb{F},k}}(X,Y)=(-1)^t \frac{R_{G_{\mathbb{F},k}}^{ev}(X)
R_{G_{\mathbb{F},k}}^{od}(Y) + R_{G_{\mathbb{F},k}}^{ev}(Y)
R_{G_{\mathbb{F},k}}^{od}(X)}
{ D^{k-\frac{1}{2}} (2i)^{t(k-3)}
<G_{\mathbb{F},k}, G_{\mathbb{F},k}>},
\end{eqnarray}
is a symmetrized sum of the product of   period polynomials  
of the normalized Hecke Eisenstein series $  G_{\mathbb{F},{k}}(\tau)$    given as followings:


\begin{prop}\label{Ei}
Take $$w_{G_{\mathbb{F}, {k}}}^{-} = \frac{\sqrt{D}\Gamma(k-1)^t   }{2^t}, \, \,
 w_{G_{\mathbb{F}, {k}}}^{+}= 
 \frac {  D^{k-\frac{3}{2}} 
 \zeta_{\mathbb{F}}(k-1)  }{(2  \pi i)^{t(k-1) }} w_{G_{\mathbb{F}, {k}}}^{-} $$
and 
\begin{eqnarray*}
 && p_{ {k}}^{+}(X )=\mathcal{N}(X)^{k-2}+(-1)^{t}, \,  \\
&& p_{ {k}}^{-}(X )=\sum_{-1 \leq n\leq k-1, n\equiv 1\pmod{2}}
 \frac{ 
\zeta_{\mathbb{F}}(1-(n+1))
 \zeta_{\mathbb{F}}(n+2-k)}{ \Gamma(n+1)^t
  \Gamma(k-n-1)^t} {\mathcal{N}(X)}^{ n}.
\end{eqnarray*}
 For $k\geq 2$
the period function  of $G_{\mathbb{F},{k}}(\tau)$ is given by
\begin{eqnarray}\label{pe}
R_{G_{\mathbb{F},{k}}}(X)=(-1)^{t }\bigl(
w_{G_{\mathbb{F}, {k}}}^{-}\cdot p_{ {k}}^- (X)+
w_{G_{\mathbb{F}, {k}}}^{+}\cdot p_{ {k}}^+ (X)
\bigr)
\end{eqnarray}
so that
$$R_{G_{\mathbb{F},{k}}}^{ev}(X)=w_{G_{\mathbb{F}, {k}}}^{+}\cdot p_{ {k}}^{+} (X) 
\mbox{ and } R_{G_{\mathbb{F},{k}}}^{od}(X) = w_{G_{\mathbb{F}, {k}}}^{-}\cdot p_{ {k}}^{-} (X).$$
\end{prop}

\begin{rmk}\label{vv}
\begin{enumerate}
\item   Like in the case of an elliptic modular form (see \cite{KZ}), the period function $R_{G_{\mathbb{F},{k}}}(X) $ is  $\frac{1}{\mathcal{N}(X)}$ times a polynomial  : 
$$R_{G_{\mathbb{F},{k}}}(X) \in \frac{1}{\mathcal{N}(X)}  \mathbb{C}[X]. $$ 

\item\label{vvv}  \cite{Geer} Note that
$\frac{D^{n-\frac{1}{2}} \zeta_{\mathbb{F}}(n)\Gamma(n)^t}
{(2\pi i)^{tn}}=\frac{\zeta_{\mathbb{F}}(1-n)}{2^t}$
and 
$ \zeta_{\mathbb{F}}(-n)=0$  for any positive even integer $n.$
\item 
\begin{enumerate}
\item \cite{Z1976}  Let 
$\mathbb{F}$ be a real quadratic field with discriminant $D.$ 
It is known that  
 $\zeta_{\mathbb{F}}(1-n)=B_{n}B_{n,\chi}$ for even positive integer $n.$
$B_r$ and $ B_{r,\chi}$ are 
 the $r$th Bernoulli number $(B_0=1, B_1=-\frac{1}{2}, B_2=\frac{1}{6}, \cdots )$ and the $r$th  twisted Bernoulli number $( B_{1,\chi}=\frac{1}{D}\sum_{a=1}^D \chi(a)a, B_{2,\chi}=\frac{1}{D}\sum_{a=1}^{D}\chi(a)a^2-\sum_{a=1}^D \chi(a) a, B_{3,\chi}= \cdots ),$  respectively.  Here 
  $\chi  \pmod{D} $ is a primitive character defined by $\chi(\cdot)=\bigl(\frac{D}{\cdot}\bigr). $  
\item { \bf{(open problem) }}
   It is well known that the  generating functions of $B_j$ and $B_{j,\chi} $    are
$$\sum_{n=0}^{\infty}B_n\frac{t^n}{n!}=\frac{te^t}{e^t-1}
 \mbox{\, and } \sum_{m\geq 0} B_{m,\chi}\frac{t^m}{m!} = \sum_{j=1}^{D}\frac{\chi(j) t e^{jt}}{e^{Dt}-1}.$$
Similarly, it will be interesting 
   to express the following    generating function
$$ \sum_{m\geq 2}^{\infty} B_m B_{m, \chi}\frac{t^m}{m!}=\sum_{m\geq 2}^{\infty}\zeta_{\mathbb{F}}(1-m)\frac{t^m}{m!},$$
as   elementary functions.
\end{enumerate}
 \end{enumerate}
\end{rmk}

\bigskip
{\pf} of Proposition \ref{Ei} :
The period polynomial of $G_{\mathbb{F},k}(\tau)$ can be computed from the definition in (\ref{Pe}) :
using $\Lambda(G_{\mathbb{F},k}, n+1)=\frac{D^{n+1}\Gamma(n+1)^t}{(2\pi)^{t(n+1)} }
 \zeta_{\mathbb{F}}(n+1)\zeta_{\mathbb{F}}(n+2-k),$  we have
\begin{eqnarray*}
&&
R_{G_{\mathbb{F}, k}}(X)
 =  \frac{(-1)^t\sqrt{D} \zeta_{\mathbb{F}}(1-k)}{2^t (k-1)^t}
 (\mathcal{N}(X)^{k-1}
+(-1)^{t}\mathcal{N}(X)^{-1})
 \\
&&  +
\frac{ (-1)^t{D}^{k-1} \Gamma(k-1)^t 
   \zeta_{ \mathbb{F}}(k-1)}{2^t (2\pi i)^{t(k-1)}}
\bigl( \mathcal{N}(X)^{k-2}+(-1)^{t}\bigr)\\
\\
&&  +
 \sum_{n=1}^{k-3}(-1)^{\frac{t(k+n-1)}{2} }
 \frac{ 
 \Gamma({k-1})^tD^{k-n-1}
}{
(2\pi)^{t(k-n-1)}  
 \Gamma(n+1)^t}
\zeta_{\mathbb{F}}(1-(n+1))
\zeta_{\mathbb{F}}(k-1-n)
\mathcal{N}(X)^{ { n}}
\end{eqnarray*}

(using the  functional equation of $\zeta_{\mathbb{F}}(s)$ of Remark \ref{vv}-(\ref{vvv}))
 \begin{eqnarray*} &&
=\frac{ (-1)^t{D}^{k-1} \Gamma(k-1)^t 
   \zeta_{ \mathbb{F}}(k-1)}{2^t(2\pi i)^{t(k-1)}}
\bigl( \mathcal{N}(X)^{k-2}+(-1)^{t}\bigr)\\
&&  +
  \frac{(-1)^t \sqrt{D} \Gamma(k-1)^t}{2^t}\sum_{-1\leq n\equiv 1\pmod{2} \leq k-1}
\frac{ \zeta_{\mathbb{F}}(1-(n+1)) 
\zeta_{ \mathbb{F}}(n+2-k)}{
\Gamma( n+1)^t \Gamma(k-n-1 )^t}
 \mathcal{N}(X)^{n}\\ &&
=
(-1)^t\bigl(\omega_{G_{\mathbb{F},k}}^{ +}p_{{k}}^+(X)+\omega_{G_{\mathbb{F}, k}}^{-}p_{{k}}^-(X)\bigr).
\end{eqnarray*}
 This completes a proof.
{\qed} 

\medskip
\subsubsection{Petersson scalar product of Eisenstein Series}\label{Eisen2}

The Petersson scalar product of a $SL_2(\mathbb{Z})$-invariant function  had been defined  by Zagier  \cite{Z-rapid} using Rankin-Selberg method.
Similarly we have 
\begin{prop}\label{innerp}
  $$ < G_{\mathbb{F},  {k}}(\tau),  G_{\mathbb{F},  {k}}(\tau) >=
 \frac{\Gamma(k-1)^t\zeta_{\mathbb{F}}(k-1) }{  (4\pi  )^{t(k-1)} }
\frac{ \zeta_{\mathbb{F}}(1-k)}{2^t}    $$
\end{prop}
{\pf} of Proposition \ref{innerp} :
Following the method   in \cite{Z-rapid}  the   Petersson norm of the Hilbert Eisenstein series  $G_{\mathbb{F}, {k}}(\tau)$ can be computed as
\begin{eqnarray*}
  < G_{\mathbb{F},  {k}}(\tau),  G_{\mathbb{F},  {k}}(\tau) >=
(-1)^{\frac{tk}{2}}(4\pi)^{-tk }\Gamma(k)^t\cdot \zeta^*_{\mathbb{F}}(k)\zeta^*_{\mathbb{F}}(2-k)
\end{eqnarray*}
where 
 $\zeta_{\mathbb{F}}^*(s):=D^{\frac{s}{2}}\pi^{-\frac{ts}{2}}\Gamma(\frac{s}{2})^t\zeta_{\mathbb{F}}(s)=\zeta_{\mathbb{F}}^*(1-s)$ (see p 57 \cite{Z1976}).
Using the   identities 
$ \Gamma(\frac{k}{2}) \Gamma(\frac{k-1}{2})= \Gamma(k-1) \sqrt{\pi} 2^{- (k-2)} $ and $\zeta_{\mathbb{F}}(k)=D^{-k+\frac{1}{2}}\frac{(2\pi i)^{tk}}{2^t\Gamma(k)^t}\zeta_{\mathbb{F}}(1-k)$ 
we get the Petersson norm of the Eisenstein series $G_{\mathbb{F},k}.$ 
{\qed}
 \\

Now  write the function  $ C_k(X,Y;\tau )$ in (\ref{kge}) as a sum of 
$ C^{cusp}_k(X,Y;\tau ) $ in  (\ref{cc})  and 
$ C^{Eis}_k(X,Y;\tau ) := R_{G_{\mathbb{F},k} }(X,Y) G_{\mathbb{F},k}(\tau) :$
 $$ C_k(X,Y;\tau ) =C^{cusp} _k(X,Y;\tau )+ C^{Eis}_k(X,Y;\tau ).$$

\begin{prop}\label{eisen} We have 
\begin{enumerate}
\item 
$C^{Eis}_k(X,Y;\tau )=(-1)^t\frac{2^t \Gamma(k-1)^t}{\zeta_{\mathbb{F}}(1-k)}
(p^+_k(X) p_k^-(Y)+p^+_k(Y) p_k^-(X)) G_{\mathbb{F},k}(\tau).$ 
\item  $C(X,Y;\tau;T)=F_{\tau}(XT,YT)F_{\tau}(T, -XYT) \mbox{
as $\tau \rightarrow (i\infty, \cdots,  i\infty)$}$ 
\end{enumerate}
\end{prop}

{\pf} of Proposition \ref{eisen} : 
\begin{enumerate}
\item   From (\ref{pee})  recall that 
$$R_{G_{\mathbb{F},k}}(X,Y)=(-1)^t \frac{R_{G_{\mathbb{F},k}}^{ev}(X)
R_{G_{\mathbb{F},k}}^{od}(Y) + R_{G_{\mathbb{F},k}}^{ev}(Y)
R_{G_{\mathbb{F},k}}^{od}(X)}
{ D^{k-\frac{1}{2}} (2i)^{t(k-3)}
<G_{\mathbb{F},k}, G_{\mathbb{F},k}>}. $$ 
So, Proposition \ref{Ei} and Proposition \ref{innerp} imply that
\begin{eqnarray*}
 &&C^{Eis}_k(X,Y;\tau )
=R_{G_{\mathbb{F},k}}(X,Y)G_{\mathbb{F},k}(\tau)\\
&& =\frac{
\omega^{+}_{G_{\mathbb{F},k}}\omega^{-}_{G_{\mathbb{F},k}}\bigl(p^+_{k}(X) p^-_{k}(Y) + p^+_{k}(Y) p^-_{k}(X) \bigr)}
{ D^{k-\frac{1}{2}} (2i)^{t(k-3)} < G_{\mathbb{F},k}, \, G_{\mathbb{F},k}>}G_{\mathbb{F},k}(\tau)
\\&&
= (-1)^t \frac{2^t \Gamma(k-1)^t}{\zeta_{\mathbb{F}}(1-k)}
(p^+_k(X) p_k^-(Y)+p^+_k(Y) p_k^-(X)) G_{\mathbb{F},k}(\tau).
\end{eqnarray*}

\item 
Using    Proposition \ref{eisen} part (1)
 the value of $C(X, Y; \tau; T)$ as
$\tau \rightarrow (i\infty, \cdots, i\infty ) $ 
 is 
\begin{eqnarray*}\label{constant} 
&&C(X,Y; (i\infty, \cdots,  i\infty);T) 
\\
&& = \frac{(\mathcal{N}(X)+\mathcal{N}(Y))(\mathcal{N}(XY)+(-1)^t)}{  \mathcal{N}(XYT)^2}+ 
\sum_{k\geq 2}C_k^{Eis}(X,Y; (i\infty, \cdots, i\infty)) 
\frac{\mathcal{N}(T)^{k-2}}{\Gamma(k-1)^t}  \nonumber \\
\\
&=&  \frac{(\mathcal{N}(X)+\mathcal{N}(Y))(\mathcal{N}(XY)+(-1)^t)}{  \mathcal{N}(XYT)^2}+ (-1)^t \sum_{k\geq 2}  
 \bigl(
p_{ {k}}^+(X)p_{ {k}}^-(Y)+p_{ {k}}^+(Y)p_{ {k}}^-(X) \bigr)
 {\mathcal{N}(T)^{k-2}}
\end{eqnarray*}
 \noindent  since 
 $G_{\mathbb{F},k}(i \infty)=\frac{\zeta_{\mathbb{F}}(1-k)}{2^t}.$ \, \, 
 On the other hand, a direct computation shows that 
\begin{eqnarray*}
&& F_{\tau}(T,-XYT) F_{\tau}(XT,YT)|_{\tau \rightarrow (i\infty, \cdots,  i\infty)}
\\
\\&&
=\frac{(\mathcal{N}(X)+\mathcal{N}(Y))(\mathcal{N}(XY)+(-1)^t)}
{   \mathcal{N}(XYT)^2}
 +(-1)^t\sum_{ k\geq 2}
\bigl( p^+_{k}(X)p^-_{k}(Y)+p^+_{k}(Y)p^-_{k}(X)\bigr)\mathcal{N}(T)^{k-2}
\end{eqnarray*}
\noindent This completes the  proof of Proposition \ref{eisen}.
\end{enumerate}
\qed

 
 \bigskip

\section{Proofs}
\subsection{{\pf} of { \bf { Theorem \ref{main}}} ({\bf{Main Theorem}})}

\noindent {\bf{ $(1)$ }}\, \, Using Theorem \ref {Rationality} with a proper choice of the Petersson norm $<f, \, f>,$  we see that
\begin{eqnarray*}
R_f(X,Y):=(-1)^t \frac{R^{ev}_f(X) R^{od}_f(Y)+R^{ev}_f(Y) R^{od}_f(X)}
{D^{k-\frac{1}{2}} \,  (2i)^{t(k-3)}  <f, \, f>} \in  \mathbb{Q}_f[X,Y], f\in S_{\mathbf{k}}.
\end{eqnarray*}

\noindent With an action of   $\sigma\in Gal(\mathbb{C}/\mathbb{Q}_f)$ by $R_{\sigma(f)}=\sigma(R_f) $ we see that
 \begin{eqnarray*}\label{ccc}
 &&  C(X,Y; \tau;T)   
 =  \frac{(\mathcal{N}(X)+\mathcal{N}(Y))(\mathcal{N}(X Y)+(-1)^t)}{  \mathcal{N}(X Y T)^2}
\\ &+& \sum_{k\geq 2}\sum_{
f\in \mathcal{B}_k
}R_f(X,Y)f(\tau)\frac{\mathcal{N}(T)^{ {k-2}}}{\Gamma(k-1)^2}\in \frac{1}{\mathcal{N}(XYT)^2}\mathbb{Q}[X,Y][[q,T]].
\end{eqnarray*}
 This proves rationality of $C(X,Y;\tau;T).$

\medskip

\noindent {\bf{$(2)$ }}\, \, To prove Theorem \ref{main} part (2)  
write the Taylor expansion (\ref{theta1}) 
\begin{eqnarray*}
 F_{\tau}(u,v) &=&  \frac{1}{
\mathcal{N}(u) }+\frac{1}{\mathcal{N}(v)} 
\\
&+& (-2)^t\sum_{k\geq 2}  \sum_{ \ell \in \mathbb{Z}_{\geq 0}^t}
\frac{ 
\mathbb{D}^{ \ell }\bigl(G_{\mathbb{F},k}(\tau)\bigr)}{(2\pi i)^{\ell}\ell ! (\ell +\mathbf{k}-\mathbf{1})!} (u^{\ell}v^{\ell+\mathbf{k}-1} +u^{\ell+\mathbf{k}-1}v^{\ell} )
 \end{eqnarray*}
or write it as
\begin{eqnarray}\label{theta3}
 F_{\tau}(u,v) =\sum_{\mathbf{h}=(h,\cdots, h), \mathbf{\ell}=(\ell_1,\cdots, \ell_2)}
g_{ {h}, \mathbf{\ell}}(\tau)(u^{\ell}v^{\ell+\mathbf{h}-\mathbf{1}}
 + u^{\ell+\mathbf{h}-\mathbf{1}}
v^{\ell})
\end{eqnarray}
with
$$g_{ {h}, \mathbf{\ell}}(\tau)=\left\{ \begin{array}{ccll}
\frac{(-2)^t}{\,  (2\pi i)^{\ell}\Gamma(\mathbf{\ell}+\mathbf{1} )  \Gamma (\mathbf{\ell}+\mathbf{h}  ) }
\mathbb{D}^{ \ell }(G_{\mathbb{F},  {h}}(\tau)), &
\mbox{\, if $  h  \geq 2, \ell_i \geq 0, i=1,2$}\\
\frac{1}{2^t}, &  \mbox{\, if $ {h}=0,    \ell  =(0, \cdots, 0)$}\\
0, & \mbox{\, otherwise}
\end{array}\right \}.$$

Next let
$$ F_{\tau}(T,  -XYT) F_{\tau}(XT,YT)=\sum_{ {k}\geq 0}
b_{\mathbf{k}}(X,Y;\tau) \frac{\mathcal{N}(T)^{ {k}- {2}}}{\Gamma(k-1)^2}.$$
 
Since   we have already checked that 
  $$C(X,Y; (i \infty, \cdots, i\infty),T)=F_{\tau }(T,  -XYT) F_{\tau}(XT,YT)|_{ \tau \rightarrow (i\infty, \cdots, i\infty)} $$ in Proposition \ref{eisen},
 it is enough to confirm  that
$$\frac{\bigl< b_{\mathbf{k}}(X,Y;\cdot ), f(\cdot)\bigr>}{\bigl< f, \, f \bigr>}
=  R_f(X,Y) \mbox {\,\, for each $ f\in \mathcal{B}^0_{k}$}. $$ 
  From the expression  (\ref{theta3}) we see

\begin{eqnarray}\label{ep}
  b_{\mathbf{k}}(X,Y;\tau)
&=&  \sum_{\tiny{
\begin{array}{ccc} \ell,\mathbf{h}, \ell', \mathbf{h}'\\ \mathbf{h}+\mathbf{h}'+ {2}(\mathbf{\ell} +\mathbf{\ell} ')=\mathbf{k}\\
\mathbf{h}=(h,\cdots, h), \mathbf{h}'=(h', \cdots, h')
\end{array}
}}g_{ {h},\mathbf{\ell}} (\tau) g_{ {h}', \mathbf{\ell}'}(\tau) \nonumber \\
&&\times
[(-XY)^{\mathbf{\ell}+\mathbf{h}-\mathbf{1}}+
(-XY)^{\mathbf{\ell}}][X^{\mathbf{\ell}'}Y^{\mathbf{\ell}'+\mathbf{h}'-\mathbf{1}}
+X^{\mathbf{\ell}'+\mathbf{h}'-\mathbf{1}}Y^{\mathbf{\ell}'}]
\end{eqnarray}

The coefficients of $\mathcal{N}(X)^{ {p}}\mathcal{N}(Y)^{ {q}}$  with $q$ or $p$ equal to $-1$ or to $k-1$ involve  only $G_{\mathbb{F},k}(\tau)$ and have already been treated in Proposition \ref{eisen}. 
Also the coefficients of $\mathcal{N}(X)^p\mathcal{N}(Y)^q$ in (\ref{ep}) is invariant under $q\leftrightarrow k-2-p$ and $q\leftrightarrow p$ so that we may assume $0\leq p<q\leq \frac{k-2}{2}.$  For such $p,q$ the coefficient of $ \mathcal{N}(X)^{ {p}}\mathcal{N}(Y)^{ {q}}$
in (\ref{ep}) equals
 \begin{eqnarray*}
&&  \sum_{\tiny{
\begin{array}{cc}
\ell,\ell'\succeq  -\mathbf{1} \\ \ell+\ell'=\mathbf{p}=(p,\cdots, p)\end{array}}
}
(-1)^{| \ell|} g_{ {k}- {p}- {q}- {1}, \mathbf{\ell}}(\tau)
   g_{ {q}- {p}+ {1},\mathbf{\ell}'}(\tau)\\
  &=&  \sum_{\tiny{
\begin{array}{cc}
\ell,\ell'\succeq  -\mathbf{1} \\ \ell+\ell'=\mathbf{p}=(p, \cdots, p)\end{array}}
}
 2^{2t}\frac{(-1)^{| \ell|}  \mathbb{D}^{\ell }(G_{\mathbb{F}, k-p-q-1}) 
 \mathbb{D}^{\ell' }(G_{\mathbb{F}, q-p+1}) }{
 \ell !  (\ell+\mathbf{k-p-q}-2)!
 {\ell'}!  (\mathbf{q}-\ell')!} 
 \\
&=&\frac{2^{2t}}{ \Gamma( {q}+\ {1})^t\Gamma( {k}-  {q}- {1})^t }[
G_{\mathbb{F},  {q}- {p}+ {1}}, G_{\mathbb{F}, {k}- {q}- {p}- {1}}]^{Hil}_{\mathbf{p}}.
\end{eqnarray*}

\medskip

 Here, $[\cdot \,\, \cdot]^{Hil}_{\mathbf{p}} $ denotes the $\mathbf{p}=(p, \cdots, p)$th Rankin-Cohen bracket (see Corollary 1 in \cite{CKR}) defined by
\begin{eqnarray*}
&\,&  [G_{\mathbb{F}, {q}-{p}+{1}}, G_{\mathbb{F},{k}-{q}-{p}-{1}}]^{Hil}_{\mathbf{p}}
 := \\
&& \frac{1}{(2\pi i)^{tp}}\sum_{\tiny{\begin{array}{ccc}
0\leq  {\ell}_i  \leq  {p}\\
 \ell=(\ell_1, \cdots, \ell_t)\\
\ell+\ell'=\mathbf{p}  \end{array} } }
 \frac{(-1)^{|\mathbf{\ell}|}\Gamma( {q}+ {1})^t
 \Gamma( {k }- {q}- {1})^t   }{  {\ell'}!
  (\mathbf{{k}+ {\ell'}-{q}- {p-2 }})!  { \ell}! (\mathbf{{q}-  {\ell}})! }
\mathbb{D}^{ \mathbf{\ell} }(G_{\mathbb{F},  {q}- {p}+ {1}})
\mathbb{D}^{ {\ell'} }(G_{\mathbb{F}, {k}- {q}- {p}- {1}}).
\end{eqnarray*}

On the other hand we recall the following result (Theorem 3 in \cite{CKR}) :

\begin{thm}\label{ckr} \cite{CKR}
 Suppose that $f(\tau)=\sum_{\mathcal{D}^{-1} \ni \nu \succ 0} a_f(\nu)e^{2\pi i tr(\nu\tau) }\in S_{\mathbf{k}}$ and $g(\tau)=\sum_{  \mathcal{D}^{-1} \ni \nu\succeq 0} a_g(\nu)e^{2\pi i tr(\nu\tau)} \in M_{\mathbf{k}_2} $ with $k=k_1+k_2+2p >2.$
Then
$$ \frac{D^{\frac{1}{2}-k_1}(2\pi i)^{t k_1}}{\Gamma(k_1)^t}
<f, [G_{\mathbb{F}, k_1}, g_{k_2}]_{\mathbf{p}}>=
 \frac{\Gamma(k-1)^t\Gamma(k_1+p)^t}
{(4\pi)^{t(k-1)}\Gamma(k_1)^t\Gamma(p+1)^t}\sum_{\nu\succ 0}
\frac{a_f(\nu)\overline{a_g(\nu)}}{\mathcal{N}(\nu)^{k-p-1}}$$
\end{thm}
\medskip

\noindent Taking $g_{k_2}(\tau)=G_{\mathbb{F},k}(\tau)$ in Theorem \ref{ckr}  we get
\begin{eqnarray*}
&&<f, [G_{\mathbb{F}, k_1}, G_{\mathbb{F}, k_2}]_{\mathbf{p}}>=
(-1)^{ \frac{tk_1}{2}} \frac{D^{k_1-\frac{1}{2}}\Gamma(k-1)^2\Gamma(k_1+p)^t}
{2^{t(k-1)}(2\pi)^{t(k+k_1-1)}\Gamma(p+1)^t}\sum_{\nu>0}
\frac{a_f(\nu)\sigma_{k_2-1}(\nu)}{\mathcal{N}(\nu)^{k-p-1}}
\\
&&
= (-1)^{ \frac{tk_1}{2}} \frac{D^{k_1-\frac{1}{2}}\Gamma(k-1)^t\Gamma(k_1+p)^t}
{2^{t(k-1)}(2\pi)^{t(k+k_1-1)}\Gamma(p+1)^t}
L(f, k-p-1) L(f, k_1+p) \\
&&
(\mbox{since $R_n(f)=i^{t(n+1)}D^{n+1}(2\pi)^{-t(n+1)} \Gamma(n+1)^t L(f,n+1)$})
\\
&&=
  \frac{  (-1)^{\frac{t(k-1)}{2}}\Gamma(k-1)^t}
{D^{k-\frac{1}{2}} 2^{t(k-1)}  \Gamma(p+1)^t \Gamma(k-p-1)^t }
R_{k-p-2}(f) R_{k_1+p-1}(f). 
\end{eqnarray*}
 And so we have  that 
\begin{eqnarray*}
 \sum_{f \in \mathcal{B}^0_k }
 \frac{R_{k-p-2}(f)  R_{k_1+ p-1}(f)  f(\tau)}{
i^ {t(k-1)} D^{k-\frac{1}{2}} 2^{t(k-3)}  <f, f>}
= \frac{   2^{2t}
\Gamma(p+1)^t
\Gamma(k-p-1)^t}
{  \Gamma(k-1)^t}
 [G_{\mathbb{F}, k_1}, 
 G_{\mathbb{F}, k_2}]^{Hil}_{\mathbf{p}}.
\end{eqnarray*}

\noindent  Now take $k_1=k-q-p-1, k_2=q+1-p $ for  $q+p\equiv 1 \pmod{2}, p, q>0, $ to get 
 
\begin{eqnarray*}
 \sum_{f\in \mathcal{B}^0_k}  (-1)^t
\frac{R_{ {k-q-2}}(f)R_{ {k-p-2}}(f) }{
  D^{k-\frac{1}{2}} (2i)^{t (k-3)}   <f, \, f>}f(\tau)
= \frac{  2^{2t} \Gamma( {k-p-1})^t  \Gamma( {p+1})^t }
{ \Gamma( {k-1})^t}
 [ G_{\mathbb{F}, {k-1-q-p}}, \, G_{\mathbb{F},  {q+1-p}}]^{Hil}_{\mathbf{p}}.
\end{eqnarray*}
Since
\begin{eqnarray*}
&& \sum_{f\in \mathcal{B}^0_k}  R_f(X,Y)f(\tau)= 
\sum_{f\in \mathcal{B}^0_k} (-1)^t 
\frac{R_f^{ev}(X) R_f^{ev}(Y)+
R_f^{ev}(Y) R_f^{ev}(X)}{ D^{k-\frac{1}{2}} (2i)^{t(k-3)}<f,\, f>  }
   f(\tau)  
\end{eqnarray*}
 we get  $$ \frac{ \bigl< \, b_{\mathbf{k}}(X,Y; \cdot), \, f(\cdot) \, \bigr>}
{ \bigl< \,f, \, f \bigl>}=
 R_f(X,Y).$$

\noindent Combining all together with Proposition \ref{eisen} we conclude that
\begin{eqnarray*}
&& C(X,Y;\tau;T)= F_{\tau}(T,-XYT)\,  F_{\tau}(XT,YT)
\end{eqnarray*}
 which completes a proof.
\qed

\medskip

\section{{\bf{Conclusion}}}
One of the main importance of   modular forms in number theory is that
spaces of modular forms are generated by those with rational Fourier coefficients.
The "period theory" gives another natural rational structure of modular forms.
A striking  result   by  Zagier  \cite{Z1991} states that   this rational information of
   modular forms
can be written as  a single product of  Kronecker series $F_{\tau}(u,v)$ which is a Jacobi form. 
The recent results in \cite{Bannai, Sp} show  that  Eisenstein-Kronecker numbers have a rich arithmetic nature,  such as a connection with the special Hecke $L$-function over imaginary quadratic fields and   Katz' two-variable p-adic Eisenstein measure. \\
In this paper,  we identified  the  Kronecker series
as  a  "Kuznetsov lifting" of holomorphic Hilbert Eisenstein series over  totally real number  fields with  strict class number 1.  This is the first case to connect    Kronecker series to  the critical values of   Hilbert modular $L$-functions over a totally real number field and it seems  worthwhile to explore the hidden arithmetic relations more.
\medskip

On the other hand, in terms of geometric interpretation, a modular form  can be  regarded as a section of a certain sheaf of differential forms on the open modular curve on a  congruence subgroup $\Gamma\subset SL_2(\mathbb{Z}).$
 By noting that the singular cohomology of the open modular curve is given by the group  cohomology $H^*(\Gamma, W)$ the comparison of de Rham and singular cohomology can give an  Eichler-isomorphism.
  Matsushima and Murakami \cite{MM} extended the results to show that the space of automorphic forms on a symmetric space $M$ is isomorphic to $H^*(M,S)$ for a certain locally constant sheaf $S$ over $M.$ The cohomology of Hilbert surfaces in terms of Hilbert cusp forms has been studied by many researchers including \cite{Geer, Yo1}. Relating the critical values of   $L$-functions of Hilbert cusp forms and cohomology   was first studied by Yoshida \cite{Yo1}. 
Following the work by Knopp \cite{Knopp} and Kohnen-Zagier \cite{KZ},  which provide us 
the  considerable new lights on Eichler-Shimura isomorphisms such as  rational structures of elliptic modular forms,  
we are able to 
 associate the space of   Hilbert modular forms  over the totally real number fields   to the parabolic cohomology group  in terms of the period polynomial by taking anti-derivative of Hilbert modular form \cite{C}.  


\bigskip

\bibliographystyle{amsplain}

\begin{thebibliography}{10}



\bibitem {BRW} \textsc{A. Babei,  L. Rolen and I.Wagner,} The Riemann hypothesis for period polynomials of Hilbert modular forms, Journal of Number Theory, Vol. 218, January 2021, page 44 - 61.


\bibitem{Bannai} \textsc{ K. Bannai and S.  Kobayashi,} 
 Algebraic theta functions and the p-adic interpolation of Eisenstein-Kronecker numbers,  Duke Mathematical  Journal,   153  (2010),  no. 2,  229 - 295. 


\bibitem{Bol} \textsc{G. Bol,}    Invarianten linearer Differentialgleichungen, Abhandlungen aus dem Mathematischen Seminar der Universit\"at Hamburg 
 16  (1949),  no.  3-4, 1 - 28. 


\bibitem{123}
\textsc{J. Brunier, G. van der Geer, G. Harder and D. Zagier,}  The $1-2-3$ of modular forms, Lectures from the Summer School on Modular Forms and their Applications held in Nordfjordeid, June 2004,  Edited by Kristian Ranestad,  Universitext. Springer-Verlag, Berlin, 2008.

\bibitem{C} \textsc{Y. Choie, } Parabolic cohomology and  Hilbert modular forms, Preprint (2021) .

 \bibitem{CKR} \textsc{Y. Choie, H. Kim and O. Richter,}    Differential operators on Hilbert modular forms, Journal of  Number Theory  122  (2007),
 no. 1, 25-36.

\bibitem{CL} \textsc{Y. Choie and M. Lee,}  Jacobi-like forms, Pseudodifferential operators, and Quasimodular forms, Springer Monographs in Mathematics  (2019).

\bibitem{CPZ}  \textsc{Y. Choie, Y. Park and D. Zagier,}   Periods of modular forms
on $\Gamma_0(N)$ and Products of Jacobi Theta functions,  Journal of the European Mathematical Society,  Vol. 21, Issue 5,   1379 - 1410 (2019).






\bibitem{Eichler}  \textsc{M. Eichler,}  Eine Verallgemeinerung der Abelschen Integrale, Mathematische Zeitschrift, vol. 67 (1957) pp 267-298.


 \bibitem{EZ} \textsc{M. Eichler and  D.  Zagier,}  The theory of Jacobi forms,  Progress in Mathematics, 55,  Birkh\"auser Boston, Inc., Boston, MA, 1985. 

\bibitem{Geer} \textsc{G. van der Geer,}   Hilbert modular surfaces,  Ergebnisse der Mathematik und ihrer Grenzgebiete (3) [Results in Mathematics and Related Areas (3)], 16, Springer-Verlag, Berlin, 1988.


\bibitem{Ha} \textsc{M. Harris,}   $L$-functions
of $2\times 2$ unitary groups and factorization of periods of Hilbert modulr forms,  Journal of American  Mathematical Society 6 (1993), no. 3, 637 - 719.


 \bibitem{Knopp} \textsc{M. Knopp,} Some new results on the Eichler cohomology of automorphic  forms, Bulletin of the American Mathematical Society,  80 (1974), 607 - 632.





\bibitem{KZ} \textsc{W. Kohnen and D. Zagier,}
Modular forms with rational periods,   Modular forms (Durham, 1983), 197-249, Ellis Horwood Series in Mathematics and its Applications: Statistics and Operational Research, Horwood, Chichester, 1984. 

\bibitem{MM} \textsc{Y. Matsushima and  S. Murakami,} 
On vector bundle valued harmonic forms and automorphic forms on symmetric Riemannian manifolds,
Annals of Mathematics (2) 78 (1963),   365 - 416. 


\bibitem{M}  \textsc{Y. Manin,}  Periods of cusp forms, and $p$-adic Hecke series,  (Russian) Matematicheskii Sbornik (N.S.)  92 (134)  (1973), 378 - 401, 503.




   \bibitem{Shimura-inv} $\blanks $, On the critical values of certain Dirichlet series and the periods of  automorphic forms, Inventiones Mathematicae
94 (1988), no. 2, 245 - 305.


 \bibitem{Shimura-Ame} $\blanks $, Algebraic relations between critical values of zeta  functions and inner products, American Journal of Mathematics,  105 (1983), no. 1, 253 - 285.


\bibitem{Shimura-Duke} $\blanks $, The special values of the zeta functions associated with Hilbert modular forms,  Duke Mathematical Journal,  45 (1978), no. 3, 637- 679.


 \bibitem{Shimura1959} $\blanks $, Sur les int\'egrales attach\'ees aux forms automorphes, Journal of the Mathematical Society of Japan   11 (1959), 291 - 311.



\bibitem{Sp} \textsc{J. Sprang, }  
Eisenstein-Kronecker series via the Poincaré bundle,
Forum of Mathematics, Sigma 7  (2019),  No.  e34, 59 pp.


\bibitem{Yo1} \textsc{H. Yoshida,} Absolute CM-periods. Mathematical Surveys and Monographs, 106,  American Mathematical Society, Providence, RI, 2003. 





\bibitem{Weil} \textsc{ A. Weil,} Elliptic functions according to Eisenstein and Kronecker, Reprint of the 1976 original,  Classics in Mathematics,  Springer-Verlag, Berlin, 1999. 


\bibitem{Z1991} \textsc{D. Zagier,} Periods of modular forms and Jacobi theta functions, Inventiones Mathematicae,  104, 449-265 (1991).

\bibitem{Z-rapid} $\blanks $, The Rankin-Selberg method for automorphic functions which are not of rapid decay, 
Journal of the Faculty of Science, University of Tokyo, Section IA, Mathematics 28 (1981), no. 3, 415 - 437 (1982).


 \bibitem{Z1976} $\blanks $,  On the values at negative integers of the Zeta-functions of a real quadratic fields, L'Enseignement Math\'ematique, 22 (1976), no. 1 - 2, 55 - 95.

\end{thebibliography}

\end{document}